# Relevant first-order logic LP<sup>#</sup> and Curry's paradox resolution

**Jaykov Foukzon**

Israel Institute of Technology, Department of Mathematics, Haifa, Israel

**Email address:**

jaykovfoukzon@list.ru

**Abstract:** In 1942 Haskell B. Curry presented what is now called Curry's paradox which can be found in a logic independently of its stand on negation. In recent years there has been a revitalised interest in non-classical solutions to the semantic paradoxes. In this article the non-classical resolution of Curry's Paradox and Shaw-Kwei's paradox without rejection any contraction postulate is proposed. In additional relevant paraconsistent logic $\check{C}_n^{\#}, 1 \leq n < \omega$, in fact, provide an effective way of circumventing triviality of da Costa's paraconsistent Set Theories $NF_n^C$

## 1. Introduction

In 1942 Haskell B. Curry presented what is now called *Curry's paradox* [1]. The paradox I have in mind can be found in a logic independently of its stand on negation. The deduction appeals to no particular principles of negation, as it is negation-free. Any deduction must use some inferential principles.

Here are the principles needed to derive the paradox.

A transitive relation of consequence: we write this by ⊢ and take ⊢ to be a relation between statements, and we require that it be transitive: if $A \vdash B$ and $B \vdash C$ then $A \vdash C$.

Conjunction and implication: we require that the conjunction operator ∧ be a greatest lower bound with respect to ⊢. That is, $A \vdash B$ and $A \vdash C$ if and only if $A \vdash B \wedge C$.

Furthermore, we require that there be a residual for conjunction: a connective → such that $A \wedge B \vdash C$ if and only if $A \vdash B \to C$.

Unrestricted Modus Ponens rule :

$$A, A \to B \vdash B \qquad (1.1)$$

Unrestricted Modus Tollens rule:

$$P \to Q, \neg Q \vdash \neg P \qquad (1.2)$$

A paradox generator: we need only a very weak paradox generator. We take the $T$ scheme in the following enthymematic form: $T[A] \wedge C \vdash A$; $A \wedge C \vdash T[A]$ for some true statement $C$. The idea is simple: $T[A]$ need not entail $A$. Take $C$ to be the conjunction of all required background constraints.

Diagonalisation. To generate the paradox we use a technique of diagonalisation to construct a statement Ψ such that Ψ is equivalent to $T[\Psi] \to A$, where $A$ is any statement you please.

Curry's paradox, is a paradox within the family of so-called paradoxes of self-reference (or paradoxes of circularity). Like the liar paradox (e.g., 'this sentence is false') and Russell's paradox, Curry's paradox challenges familiar naive theories, including naïve truth theory (unrestricted $T$-schema) and naive set theory (unrestricted axiom of abstraction), respectively. If one accepts naive truth theory (or naive set theory), then Curry's paradox becomes a direct challenge to one's theory of logical implication or entailment. Unlike the liar and Russell paradoxes Curry's paradox is negation-free; it may be generated irrespective of one's theory of negation.

There are basically two different versions of Curry's paradox, a truth-theoretic (or proof-theoretic) and a set-theoretic version; these versions will be presented below.

*Truth-theoretic version.*

Assume that our truth predicate satisfies the following $T$-schema:

$$T[A] \leftrightarrow A$$

Assume, too, that we have the principle called Assertion (also known as pseudo modus ponens):

Assertion: $(A \wedge (A \to B)) \to B$

By diagonalization, self-reference we can get a sentence $C$ such that $C \leftrightarrow (T[C] \to F)$, where $F$ is anything you like. (For effect, though, make $F$ something obviously false, e.g. $F \equiv 0 = 1$) By an instance of the $T$-schema ("$T[C] \leftrightarrow C$"), we immediately get: $T[C] \leftrightarrow (T[C] \to F)$,

Again, using the same instance of the $T$-Schema, we can substitute $C[T,F]$ for $T[C]$ in the above to get (1).

1  ⊢ $C[T,F] \leftrightarrow (C[T,F] \to F)$  [by $T$-schema and Substitution]
2  ⊢ $(C[T,F] \wedge (C[T,F] \to F)) \to F$ [by Assertion]
3  ⊢ $(C[T,F] \wedge C[T,F]) \to F$ [by Substitution, from 2]
4  ⊢ C[T,F] → F [by Equivalence of $C$ and $C \wedge C$, from 3]
5  ⊢ $C[T,F]$ [by Modus Ponens, from 1 and 4]
6  ⊢ $F$ [by Modus Ponens, from 4 and 5]

Letting $F$ be anything entailing triviality Curry's paradox quickly 'shows' that the world is trivial.

*Set-Theoretic Version*

The same result ensues within naive set theory. Assume, in particular, the

(unrestricted) axiom of abstraction (or naive comprehension (NC)):

Unrestricted Abstraction: $x \in \{x|A(x)\} \leftrightarrow A(x)$.

Moreover, assume that our conditional, $\to$, satisfies Contraction (as above),

which permits the deduction of $(s \in s \to A)$ from

$$s \in s \to (s \in s \to A).$$

In the set-theoretic case, let $C[F], \{x|x \in x \to F\}$, where $F$ remains as you please (but something obviously false, e.g. $F \equiv 0 = 1$). From here we reason thus:

(1) ⊢ $x \in C[F] \leftrightarrow (x \in x \to F)$ [by Unrestricted Abstraction]
(2) ⊢ $C[F] \in C[F] \leftrightarrow (C[F] \in C[F] \to F)$ [by Universal Specification, from 1]
(3) ⊢ $C[F] \in C[F] \to (C[F] \in C[F] \to F)$ [by Simplification, from 2]
(4) ⊢ $C[F] \in C[F] \to F$ [by Contraction, from 3]
(5) ⊢ $C[F] \in C[F]$ [by Unrestricted Modus Ponens, from 2 and 4]
(6) ⊢ $F$ [by Unrestricted Modus Ponens, from 4 and 5]

So, coupling Contraction with the naive abstraction schema yields, via Curry's paradox, triviality.

This is a problem. Our true $C[F]$ entails an arbitrary $F$. This inference arises independently of any treatment of negation. The form of the inference is reasonably well known. It is Curry's paradox, and it causes a great deal of trouble to any non-classical approach to the paradoxes. In the next sections we show how the tools for Curry's paradox are closer to hand than you might think.

## 2. Relevant First-Order Logics in General

Relevance logics are non-classical logics [2]-[15]. Called "relevant logics" in Britain and Australasia, these systems developed as attempts to avoid the paradoxes of material and strict implication. It is well known that relevant logic does not accept an axiom scheme $A \to (\neg A \to B)$ and the rule $A, \neg A \vdash B$. Hence, in a natural way it might be used as basis for contradictory but non-trivial theories, i.e. paraconsistent ones. Among the paradoxes of material implication are: $p \to (q \to p), \neg p \to (p \to q), (p \to q) \vee (q \to r)$. Among the paradoxes of strict implication are the following: $(p \wedge \neg p) \to q, p \to (q \to q), p \to (q \wedge \neg q)$. Relevant logicians point out that what is wrong with some of the paradoxes (and fallacies) is that is that the antecedents and consequents (or premises and conclusions) are on completely different topics. The notion of a topic, however, would seem not to be something that a logician should be interested in — it has to do with the content, not the form, of a sentence or inference. But there is a formal principle that relevant logicians apply to force theorems and inferences to "stay on topic". This is the variable sharing principle. The variable sharing principle says that no formula of the form $A \to B$ can be proven in a relevance logic if $A$ and $B$ do not have at least one propositional variable (sometimes called a proposition letter) in common and that no inference can be shown valid if the premises and conclusion do not share at least one propositional variable.

## 3. Curry's Paradox Resolution Using Canonical Systems of Relevant Logic

In the work of Anderson and Belnap [3] the central systems of relevance logic were the logic E of relevant entailment and the system R of relevant implication. The relationship between the two systems is that the entailment connective of E was supposed to be a strict (i.e. necessitated) relevant implication. To compare the two, Meyer added a necessity operator to R (to produce the logic NR).

It well known in set theories based on strong relevant logics, like E and R, as well as in classical set theory, if we add the naive comprehension axiom, we are able to derive any formula at all. Thus, naive set theories based on systems such as E and R are said to be "trivial" by Curry Paradox.

The existence of this paradox has led Grishen, Brady, Restall, Priest, and others to abandon the axiom of contraction which we have dubbed

K: $((A \to (A \to B)) \to (A \to B))$.

Brady has shown that by removing contraction, plus some other key theses, from R we obtain a logic that can accept naive comprehension without becoming trivial [4],[16],[17].

However, it is not just K that we must avoid. Shaw-Kwei [21] shows that a variant of Curry's paradox can trivialise a chain of weaker naive truth theories. Let us use the notations

$$\varphi \to_{(0)} \psi \text{ and } \varphi \to_{(n+1)} \psi$$

to mean $\psi$ and $\left(\varphi \to (\varphi \to_{(n)} \psi)\right)$ correspondingly.

Then the following axioms also lead to triviality

$$K_n: \left(\varphi \to (\varphi \to_{(n)} \psi)\right) \to \left(\varphi \to_{(n)} \psi\right).$$

We choose now a sentence $\gamma_n$ via the diagonal lemma, that satisfies [22]:

$$\gamma_n \leftrightarrow \left(Tr(\widehat{\gamma_n}) \to_{(n)} \varphi\right),$$

where the notations $\widehat{\circ}$ to mean an fixed Godel numbering.

Then by full intersubstitutivity one obtain the equivalence

$$E_n: \gamma_n \leftrightarrow \left(\gamma_n \to_{(n)} \varphi\right),$$

which by postulate $K_n$ reduces to $\left(\gamma_n \to_{(n)} \varphi\right)$ and by $E_n$ to $\gamma_n$. But from $\gamma_n$ and $\gamma_n \to_{(n)} \varphi$ one can deduce $\varphi$ by n applications of unrestricted modus ponens (1.1). For example, a natural implicational logic without contraction is Łukasiewicz's 3-valued logic: Ł$_3$. Although logic Ł$_3$ does not contain K, it does contain $K_2$. In general the n+1-valued version of Łukasiewicz logic, Ł$_{n+1}$, validates $K_n$ and is thus unsuitable for the same reason [22],[23].

However, it well known that contraction is not the only route to triviality. There are logics which are contraction free that still trivialize naive comprehension schema (NC) [18]. Abelian logic with axiom of relativity which we have dubbed

$$R: \left((p \to q) \to q\right) \to p.$$

Let $a = \{x | \varphi(x)\}$ and $\varphi(x) = p \to x \in x$. Then as instance of NC one obtain $(p \to a \in a) \to a \in a$. Thus we obtain
  (1) $\vdash (p \to a \in a) \to a \in a$ [by NC]
  (2) $\vdash \left((p \to a \in a) \to (a \in a)\right) \to p$ [by instance of R]
  (3) $\vdash p$ [by 1,2 and Unrestricted Modus Ponens (1.1)].

## 4. Relevant First-Order Logic LP$^\#$

In order to avoid the results mentioned in II and III, one could think of restrictions in initial formulation of the rule Unrestricted Modus Ponens (1.1). The postulates (or their axioms schemata) of propositional logic LP$^\#$[V] are the following [19]:
  I. Logical postulates:
  (1) $A \to (B \to A)$,
  (2) $(A \to B) \to ((A \to (B \to C)) \to (A \to C))$,
  (3) $A \to (B \to A \wedge B)$,
  (4) $A \wedge B \to A$,
  (5) $A \wedge B \to B$,
  (6) $A \to (A \vee B)$,
  (7) $B \to (A \vee B)$,
  (8) $(A \to C) \to ((B \to C) \to (A \vee B \to C))$,
  (9) $A \vee \neg A$,
  (10) $B \to (\neg B \to A)$.
  II. Restricted Modus Ponens rule:

$$A, A \to B \vdash B \text{ iff } A \notin V \qquad (1.3)$$

or

$$A, A \to B \vdash B \text{ iff } B \notin V \qquad (1.4)$$

which we have write for short

$$A, A \to B \vdash_r B \text{ or } A, A \to B \vdash_{r,V} B.$$

## 5. Curry's Paradox and Shaw-Kwei's Paradox Resolution Using Relevant First-Order Logic LP$^\#$

In my paper [19] was shown that by removing only Unrestricted Modus Ponens rule (1.1) (without removing contraction etc.), plus some other key theses, from classical logic we obtain a logic that can accept naive comprehension without becoming trivial.

Let us consider Curry's paradox in a set theoretic version using Relevant First-Order Logic LP$^\#$ with Restricted Modus Ponens rule (1.3). Let $C[F] = \{x | x \in x \to F\}$ and $\alpha[F]$ is a closed a well formed formula of ZFC (cwff) such that: $\alpha[F] \leftrightarrow C[F] \in C[F]$. We assume now $Con(ZFC)$ and denote by $\Delta$ a set of all cwff such that $\beta \in \Delta \leftrightarrow \neg Con(ZFC + \beta)$. Let us denote by symbol $W_\Delta$ a set

$$W_\Delta = \{ C[F] | F \in \Delta \}$$

We set now in (1.3). $V = W_\Delta$ From definition above we obtain the Restricted Modus Ponens rule:

$$A, A \to B \vdash B \text{ iff } A \notin W_\Delta. \qquad (1.5)$$

Let $F \in \Delta$. From here we reason thus:
  (1) $\vdash x \in C[F] \leftrightarrow (x \in x \to F)$ [by Unrestricted Abstraction]
  (2) $\vdash C[F] \in C[F] \leftrightarrow (C[F] \in C[F] \to F)$ [by Universal Specification, from 1]
  (3) $\vdash C[F] \in C[F] \to (C[F] \in C[F] \to F)$ [by Simplification, from 2]
  (4) $\vdash C[F] \in C[F] \to F$ [by Contraction, from 3]
  (5) $\nvdash_r C[F] \in C[F]$ [by Restricted Modus Ponens (1.5), from 2 and 4]

Let us denote by symbol $\widetilde{W}_\Delta$ a set

$$\widetilde{W}_\Delta = \{ C[F] | F \notin \Delta \}.$$

Therefore

$$A, A \to B \vdash B \text{ iff } A \in \widetilde{W}_\Delta \qquad (1.6)$$

Let $F \notin \Delta$. From here we reason thus:
  (1) $\vdash x \in C[F] \leftrightarrow (x \in x \to F)$ [by Unrestricted Abstraction]
  (2) $\vdash C[F] \in C[F] \leftrightarrow (C[F] \in C[F] \to F)$ [by Universal Specification, from 1]
  (3) $\vdash C[F] \in C[F] \to (C[F] \in C[F] \to F)$ [by Simplification, from 2]
  (4) $\vdash C[F] \in C[F] \to F$ [by Contraction, from 3]
  (5) $\vdash C[F] \in C[F]$ [by Restricted Modus Ponens (1.6), from 2 and 4]
  (6) $\vdash F$ [by Restricted Modus Ponens (1.6), from 4 and 5]

Let us consider now Curry's paradox in a set theoretic version using Relevant First-Order Logic LP$^\#$ with Restricted Modus Ponens rule (1.4). We set now in (1.4). $V = W_\Delta$ From definition above we obtain the Restricted Modus Ponens rule:

$$A, A \to B \vdash B \text{ iff } B \notin W_\Delta \quad (1.7)$$

Let $F \in \Delta$. From here we reason thus:
(1)   $\vdash x \in C[F] \leftrightarrow (x \in x \to F)$ [by Unrestricted Abstraction]
(2)   $\vdash C[F] \in C[F] \leftrightarrow (C[F] \in C[F] \to F)$ [by Universal Specification, from 1]
(3)   $\vdash C[F] \in C[F] \to (C[F] \in C[F] \to F)$ [by Simplification, from 2]
(4)   $\vdash C[F] \in C[F] \to F$ [by Contraction, from 3]
(5)   $\nvdash_r C[F] \in C[F]$ [by Restricted Modus Ponens (1.7), from 2 and 4]

Let us consider now Curry's paradox in a set theoretic version using Abelian logic with axiom of relativity and Restricted Modus Ponens (1.4). We set now in (1.4). $V = \Delta$ From definition above we obtain the Restricted Modus Ponens rule:

$$A, A \to B \vdash B \text{ iff } B \notin \Delta \quad (1.8)$$

Let $C[F] = \{x | \varphi(x)\}$ and $\varphi(x) = F \to x \in x$ and let $F \in \Delta$. Then as instance of NC one obtain $(F \to C[F] \in C[F]) \to C[F] \in C[F]$. Thus we obtain
(1)   $\vdash (F \to (C[F] \in C[F])) \to (C[F] \in C[F])$ [by NC]
(2)   $\vdash ((F \to (C[F] \in C[F])) \to (C[F] \in C[F])) \to F$ [by instance of $R$]
(3)   $\nvdash_r F$ [by 1,2 and Restricted Modus Ponens (1.7)].

Let us consider now Curry's paradox in a truth theoretic version using Relevant First-Order Logic LP# with Restricted Modus Ponens rule (1.4). We set now in (1.4). $V = \Delta$ From definition above we obtain the Restricted Modus Ponens rule:

$$A, A \to B \vdash B \text{ iff } B \notin \Delta. \quad (1.9)$$

By diagonalization, self-reference we can get a sentence $C$ such that $C \leftrightarrow (T[C] \to F)$, where $F \in \Delta$.

By an instance of the $T$-schema ("$T[C] \leftrightarrow C$") we immediately get: $T[C] \leftrightarrow (T[C] \to F)$,

Again, using the same instance of the $T$-Schema, we can substitute $C[T, F]$ for $T[C]$ in the above to get (1).
(1)   $\vdash C[T, F] \leftrightarrow (C[T, F] \to F)$ [by $T$-schema and Substitution]
(2)   $\vdash (C[T, F] \wedge (C[T, F] \to F)) \to F$ [by Assertion]
(3)   $\vdash (C[T, F] \wedge C[T, F]) \to F$ [by Substitution, from 2]
(4)   $\vdash C[T, F] \to F$ [by Equivalence of $C$ and $C \wedge C$, from 3]
(5)   $\vdash C[T, F]$ [by Restricted Modus Ponens (1.9), from 1 and 4]
(6)   $\nvdash_r F$ [by Restricted Modus Ponens (1.9), from 4 and 5].

It easy to see that by using logic with appropriate restricted modus ponens rule (1.4) Shaw-Kwei's paradox disappears by the same reason.

# 6. The Resolution of ω-Inconsistency Problem for the Infinite Valued Łukasiewicz Logic Ł∞. Logic $LP^{\#}_{\{\omega\}}$.

It well known that in the infinite valued Łukasiewicz logic, Ł∞, every instance of $K_n$ is invalid, and in fact Ł∞ can consistently support a naive truth predicate [23]-[24]. However, Ł∞ is plagued with an apparently distinct problem – it is $\omega$-inconsistent. This fact was first shown model theoretically by Restall in [25] and demonstrated a proof theoretically by Bacon in [24].

An classical extension of Peano Arithmetic is said to be $\omega$–inconsistent iff

$$\vdash \varphi[n/x] \text{ for each n, but } \vdash \exists x \neg \varphi[x] \quad (1.10)$$

Note that while an $\omega$-inconsistent theory is not formally inconsistent. However $\omega$-inconsistency is generally considered to be an undesirable property, generally considered to be an undesirable property. It is generally considered undesirable if the theory becomes inconsistent in $\omega$–logic. In other words, if it cannot be consistently maintained in the presence of the infinitary $\omega$-rule:

$$\{\varphi[n/x] | n \in \omega\} \vdash \forall x \, \varphi[x] \quad (1.11)$$

Clearly $\omega$-inconsistency entails inconsistency with the $\omega$–rule (1.11), but the converse does not hold in general. We have dubbed any Logic LP# with the $\omega$–rule (1.11) by $LP^{\#}_{\{\omega\}}$.

Definition 6.1. [23]. Weak $\omega$–inconsistency means:

$$\varphi[n/x] \vdash \text{ for each } n, \text{ but } \vdash \exists x \varphi[x]. \quad (1.12)$$

Definition 6.2. [23]. Strong $\omega$–inconsistency means:

$$\vdash \varphi[n/x] \text{ for each } n, \text{ but } \vdash \exists x (\varphi[x] \to \bot). \quad (1.13)$$

Note that without the rule of reduction one cannot derive strong $\omega$-inconsistency from weak $\omega$-inconsistency [23].

Definition 6.3. [23]. By a classical "naive truth theory" (CNTT) we shall mean any set of first order sentences in the language of arithmetic with a truth predicate which, in addition to being closed under modus ponens, has the following properties:
(1) Standard syntax: it contains all the arithmetical consequences of classical Peano arithmetic.
(2) Intersubstitutivity: it contains $\varphi$ if and only if it contains $\varphi \left[ Tr\left(\widetilde{\psi}\right)/\psi \right]$ for any sentence $\varphi$.
(3) Compositionality: it contains $Tr(x) \to Tr(y)$ if and only if it contains $Tr(x \to y)$.
(4) Unrestricted Modus Ponens rule: it closed under unrestricted modus ponens rule (1.1).

If $\varphi \vdash \psi$ then $\exists x \varphi \vdash \exists x \psi$, $(\varphi \to \exists x \psi) \vdash \exists x (\varphi \to \psi)$.

Note that by using the diagonal lemma we can construct a sentence $\gamma$ satisfying

$$\gamma \leftrightarrow \exists n Tr\left(f(n, \widetilde{\gamma})\right), \quad (1.14)$$

where the notations $\widetilde{\circ}$ to mean an fixed Gödel numbering and a function $f$ is defined arithmetically by recursion [23]:

$f(0,x) = x \dot{\rightarrow} \bot$ and $f(n+1,x) = x \dot{\rightarrow} f(n,x)$.

**Theorem 6.1.** [23]. Any classical naive truth theory closed under (1),(2),(3), (i) and (ii) can prove $\gamma$.

**Theorem 6.2.** [23]. Any naive truth theory closed under (i) and (ii) is weakly $\omega$-inconsistent.

Proof. By theorem 6.1 one obtain

$$\text{CNTT} \vdash \exists n Tr\left(f(n,\widehat{\gamma})\right).$$

By arithmetic and full intersubstitutivity we obtain that

$$Tr\left(f(n,\widehat{\gamma})\right) \vdash \gamma \rightarrow_{(n)} \bot.$$

Since we have $\vdash \gamma$ by theorem 2.1, by n applications of unrestricted modus ponens we obtain

$$\gamma \rightarrow_{(n)} \bot \vdash \bot. \quad (1.15)$$

So we have in general $Tr\left(f(n,\widehat{\gamma})\right) \vdash$ for any n, and

$$\vdash \exists n Tr\left(f(n,\widehat{\gamma})\right).$$

**Theorem 6.3.** [23]. Any naive truth theory closed under (i) and (ii) is strongly $\omega$-inconsistent.

**Theorem 6.4.** [25]. Infinitely valued Łukasiewicz logic, Ł$_\infty$, is strongly $\omega$-inconsistent.

**Definition 6.4.** By a non-classical or generalized "naive truth theory" (GCNTT) we shall mean any set of first order sentences in the language of arithmetic with a truth predicate which, in addition to being closed under modus ponens, has the following properties:

Standard syntax: it contains all the arithmetical consequences of classical Peano arithmetic.

Intersubstitutivity: it contains $\varphi$ if and only if it contains $\varphi\left[Tr\left(\widehat{\psi}\right)/\psi\right]$ for any sentence $\varphi$.

Compositionality: it contains $Tr(x) \rightarrow Tr(y)$ if and only if it contains $Tr(x \dot{\rightarrow} y)$.

Infinitary $\omega$-rule: $\{\varphi[n/x] | n \in \omega\} \vdash \forall x\, \varphi[x]$.

Restricted Modus Ponens rules: it closed under restricted modus ponens rule (1.3) or (1.4)

**Definition 6.5.** Weak $\omega$-consistency means:

$$\vdash \varphi[n/x] \text{ for each } n, \text{ but } \nvdash \exists x(\varphi[x] \rightarrow \bot) \quad (1.16)$$

**Definition 6.6.** [23]. Strong $\omega$-consistency means:

$$\vdash \varphi[n/x] \text{ for each } n, \text{ but } \nvdash \exists x \neg(\varphi[x]) \quad (1.17)$$

**Theorem 6.5.** Any consistent GCNTT closed under restricted modus ponens rule (1.3) is strongly $\omega$-consistent.

# 7. Applications to da Costa's Paraconsistent Set Theories.

da Costa [27] introduced a Family of paraconsistent logics $C_n, 1 \leq n \leq \omega$, with unrestricted modus ponens rule (1.3) [28], designed to be able to support set theories $NF_n^C$, respectively, $1 \leq n \leq \omega$, incorporating unrestricted **Comprehension Schema:**

$$\vdash \exists x \forall y[x \in y \leftrightarrow F(x)], \quad (1.18)$$

where $F(x)$ is any formula in which $y$ is not free but $x$ may be, and $F(x)$ does not contain any sub formula of the form $A \rightarrow B$

**Axiom of Extensionality:**

$$\vdash \forall x \forall y[\forall z(z \in x \leftrightarrow z \in y) \rightarrow x = y] \quad (1.19)$$

Since Russell's paradox could be reproduced in these set theories, their underlying logics in the absence classical rule $A, \neg A \nvdash B$ had to be capable of tolerating such theorems as $\vdash \mathcal{R} \in \mathcal{R} \leftrightarrow \neg \mathcal{R} \in \mathcal{R}$ without collapse into triviality [29] but which is hardly less disastrous $\vdash \forall x \forall y[(x \in y) \land (x = y)]$.

**Definition 7.1.** [29]. $\sim A$ iff $A \rightarrow \forall x \forall y[(x \in y) \land (x = y)]$.

**Theorem 7.1.** [29] In $NF_n^C$, negation $\sim$ is a minimal intuitionistic negation.

**Theorem 7.2.** [29]. (Cantor's Theorem)
$\vdash \forall \alpha [\sim (\alpha = P(\alpha))]$.

**Definition 7.2.** [29]. The universal set $V$ is defined as: $\forall x[x \in V \leftrightarrow (x = x)]$

**Theorem 7.3.** [29]. (Cantor's Paradox)
$\vdash (V = P(V)) \land \sim (V = P(V))$.

**Theorem 7.4.** [29].

(i) $\quad \forall x \forall y[(x = y) \land \sim(x = y)]$,

(ii) $\quad \forall x \forall y[(x \in y) \land \sim(x \in y)]$,

(iii) $\quad \forall x \forall y[(x \in x) \land \sim(x \in x)]$.

Proof. (i). By theorem 7.3 one obtain

$$\vdash (V = V) \land \sim(V = V). \quad (1.20)$$

From (1.20) and definition 7.1 one obtain

$$V = V \rightarrow \forall x \forall y[(x \in y) \land (x = y)]. \quad (1.21)$$

Therefore, as $V = V$, then $\vdash \forall x \forall y(x = y)$ and $\vdash \forall x \forall y(x \in y)$.

Note that statement (i) of the theorem 7.4 is called *paradox of identity*.

**Definition 7.3.** Let us define paraconsistent da Costa type logics $C_n^\#, 1 \leq n < \omega$, with restricted modus ponens rule such that

$$A, A \rightarrow B \vdash B \text{ iff } B \notin V, \quad (1.22)$$

$$[\forall x \forall y(x = y)] \land [\forall x \forall y(x \in y)] \in V \quad (1.23)$$

for support set theories $NF_n^\#$, respectively, $1 \leq n < \omega$, incorporating unrestricted Comprehension Schema (1.18).

From the proof of the theorem 7.3 it follows directly that logics $C_n^\#, 1 \leq n < \omega$, in fact, provide an effective way of circumventing paradox of identity.

Arruda in [29] introduced a Family of set theories $ZF_n$, $1 \leq n \leq \omega$, in which any canonical axiom of $ZFC$: the axiom of pairing, axiom of union etc., are postulated in general and in which also postulated the existence of the Russell's set $\mathcal{R}$.

Definition 7.4. [29]. $\neg_n^* A$ iff $\neg A \wedge A^{(n)}$.

Note that $\neg_n^*$ is a classical negation.

Theorem 7.5. [29]. Any set theories $ZF_n$, $1 \leq n \leq \omega$ are trivial.

Proof. By axiom of separation there exist subset $\mathcal{R}_n$ of $\mathcal{R}$ such that

(1) $\forall x [x \in \mathcal{R}_n \leftrightarrow (x \in \mathcal{R}) \wedge (x \in x)^{(n)}]$. From (1) we obtain

(2) $\mathcal{R}_n \in \mathcal{R}_n \leftrightarrow \neg(\mathcal{R}_n \in \mathcal{R}_n) \wedge (\mathcal{R}_n \in \mathcal{R}_n)^{(n)}$.

From (2) we obtain

$$(\mathcal{R}_n \in \mathcal{R}_n) \wedge \neg_n^*(\mathcal{R}_n \in \mathcal{R}_n). \quad (1.24)$$

But formula (1.24) trivializes the system $ZF_n$.

Definition 7.5. Let us define paraconsistent da Costa type logics $\check{C}_n^\#, 1 \leq n < \omega$, with restricted modus ponens rule such that

$$A, A \to B \vdash B \text{ iff } A \neg_n^* A \notin V_n. \quad (1.25)$$

$$(\mathcal{R}_n \in \mathcal{R}_n) \wedge \neg_n^*(\mathcal{R}_n \in \mathcal{R}_n) \in V_n \quad (1.26)$$

for support set theories $\widetilde{ZF}_n^\#$, respectively, $1 \leq n < \omega$, incorporating unrestricted Comprehension Schema (1.18).

From the proof of the theorem 7.5 it follows directly that logic $\check{C}_n^\#, 1 \leq n < \omega$, in fact, provide an effective way of circumventing Russell's paradox.

Arruda and da Costa [30] introduced a Family of sentential logics, $J_1$ to $J_5$, designed to be able to support set theories, respectively $ZF_1$ to $ZF_5$, incorporating an unrestricted Comprehension Schema (1.18). These $J$ logics are interesting in that they do not have modus ponens, but still seem to contain a lot of theorems that might be expected if modus ponens was included.

Theorem 7.6. [32]. $\vdash A$ is a theorem of positive intuitionistic logic if and only if $\to A$ is a theorem of $J_1$.

The basic version of Curry's paradox shows that any such set theory is trivial if its underlying logic contains the rules of Unrestricted Modus Ponens (1.1) and Contraction, in addition to the usual Instantiation rules for the quantifiers and Simplification. Arruda and da Costa instead constructed their $J$-systems without *modus ponens*. Arruda and da Costa [27] announced that $A \equiv \neg A \vdash B \supset C$ is derivable in $J_2$ to $J_5$ for all formulas $A$, $B$ and $C$. Consequently, by Russell's paradox, the set theories: $ZF_1$ to $ZF_5$ contain $\vdash B \supset C$ for all $B$ and $C$. In the absence of modus ponens, this does not quite amount to triviality. It is rather a variant which can be called $\supset$-triviality, but which is hardly less disastrous: $\forall x \forall y(x = y)$ directly follows by Axiom of Extensionality (1.19). Noting only that $A \equiv \neg A \vdash B \supset C$ is not similarly derivable in $J_1$. Arruda and da Costa [31] left open the question whether the sole remaining set theory $ZF_1$ is acceptably non-trivial, and thus whether the strategy of restricting *modus ponens* in the manner of the $J$-systems does in fact provide an effective way of circumventing Curry's paradox. These questions answered in the negative by the following variant of the Russell's paradox [33]:

$$\mathcal{R} \in \mathcal{R} \equiv (\mathcal{R} \in \mathcal{R} \supset C). \quad (1.27)$$

Theorem 7.7. [33]. $ZF_1$ is $\supset$-trivial.

In addition to Contraction, Simplification and Instantiation rules, $J_1$ contains the rules of Weakening, $B \vdash A \supset B$, and Transitivity, $A \supset B, B \supset C \vdash A \supset C$.

Definition 7.6. Let us define paraconsistent logic $J_1^\#$, with restricted Weakening rule such that

$$B \vdash A \supset B \text{ iff } B \notin V, \quad (1.28)$$

$$(\mathcal{R} \in \mathcal{R} \supset C) \in V. \quad (1.29)$$

For support set theory $ZF_1^\#$, incorporating unrestricted Comprehension Schema (1.18). From the proof of the theorem 7.7 it follows directly that logic $J_1^\#$ in fact, provide an effective way of circumventing Curry's paradox.

## 8. Conclusions

We pointed out that appropriate resolution of Curry's Paradox and Shaw-Kwei's paradox resolution can be given without rejection any contraction postulate. In additional logic $\check{C}_n^\#, 1 \leq n < \omega$, in fact, provide an effective way of circumventing triviality of da Costa's Paraconsistent Set Theories.